# Best one-sided algebraic approximation by average modulus

Raheam A. Al-Saphory[1,*], Abdullah A. Al-Hayani[2] and **Alaa A. Auad[3]

[1,2] Department of Mathematics; College of Education for Pure Sciences;
Tikrit University, Salahaddin; IRAQ.

[3] Department of Mathematic; College of Education for Pure Sciences University of Anbar
; Ramadi; IRAQ.

*E-mail: *saphory@tu.edu.iq          ** alaa.adnan.auad@uoanbar.edu.iq

# Best one-sided algebraic approximation by average modulus

**Raheam A. Al-Saphory[1,*], Abdullah A. Al-Hayani[2] and **Alaa A. Auad[3]**

[1, 2] Department of Mathematics; College of Education for Pure Sciences; Tikrit University, Salahaddin; IRAQ.
[3]Department of Mathematic; College of Education for Pure Sciences University of Anbar; Ramadi; IRAQ.

*E-mail: *saphory@tu.edu.iq          ** alaa.adnan.auad@uoanbar.edu.iq

**Abstract**

The aim of this work is to introduced the concept of the best one-sided approximation of unbounded functions in weighted space by using algebraic operators in terms the average modulus of smoothness. We also show an estimate of the degree of best one-sided approximation of unbounded functions in the terms of average modulus of smoothness.

**Keywords:** weighted space, algebraic polynomial unbounded function, and average modulus of smoothness.

## 1. INTRODUCTION

Ronald [1] in [1968] studied the problem of approximation in a normed linear space has associated with it a dual problem of maximizing functionals. Thus, Doronin and Ligun [2] discussed the problem of the one-sided approximation of functions by n-dimensional subspaces and find the exact value of the best one-sided approximation of the class $W^rL_1$. So, Nurnberger [3] studied Unity in one-sided L1-approximation and quadrature formulae. Babenko and Glushko [4] studied the problem of the uniqueness of elements of the best approximations in the space $L_1[a,b]$ and expressed the problem of the best approximation and the best $(\alpha, \beta)$-approximation of continuous functions and the problem of the best one-sided approximation of continuously differentiable functions. Thus, Yang [5] presented one sided $L_p$ norm and best approximation in one sided $L_p$ norm. So, Motornyi and Sedunova [6] found the best one-sided approximations of the class $W^1_\infty$ of differentiable functions by algebraic polynomials in $L_1$ space. Thus, Bustamante et al. [7] studied polynomials of the best one-sided approximation to a step function on [-1,1] and they proved that polynomials are obtained by Hermite interpolation at the zeros of some quasi orthogonal Jacobi polynomial. So, Alexander [8] in 2016 found the Best One-Sided Approximation of Some Classes of Functions of Several Variables by Haar Polynomials defined using modulus of continuity $\omega(f, t)$ and $\omega_{\rho i}(f, \delta)$ and in the same year, Alaa and Mousa [9] studied some positive factors for the one-sided approximation of the infinite functions in the weighted space $L_{p,\alpha}(X)$ and give an estimate of the degree of the best one-sided approximation in terms of the mean continuity coefficient Thus, Sedunova [10] studied the best one-sided approximation for the class of differentiable functions by algebraic polynomials in the mean. Also, Jianbo and Jialin [11] presened the Strong and Weak Convergence Rates of a Spatial Approximation for Stochastic Partial Differential Equation with One-sided Lipschitz Coefficient. Furthermore, Raad et al. [12] found the best one-sided multiplier approximation of unbounded functions by trigonometric polynomials in term of averaged modulus.

## 2. PRILIMINARIES

Continuing our previous investigations on polynomial operators for one-sided approximation to unbounded functions in weighted space (see [5]), it is the aim of this paper to develop a notion of direct estimation polynomial approximation with constructs which fits, to gather with results (see [8] and [9]) for unbounded function approximation processes.

Let $X = [0,1]$, we denoted by $L_{p,\beta}(X), 1 \leq p < \infty$ be the space of all unbounded functions that defined on X, with every function in this space has the norm given by

$$\|\rho\|_{p,\beta} = \left(\int_X |\rho(x)|^p\right)^{\frac{1}{p}} < \infty. \qquad (1)$$

Let W be the suitable set of all weight functions on X, such that $|\rho(x)| \leq M\beta(x)$, where $M$ is positive real number and

$\beta : X \to \mathbb{R}$ weight function, which are equipped with the following norm

$$\|\rho\|_{p,\beta} = \left(\int_X \left|\frac{\rho(x)}{\beta(x)}\right|^p\right)^{\frac{1}{p}} < \infty. \qquad (2)$$

The local modulus of continuity of a function $\rho : [0,1] \to \mathbb{R}$ in the pint $x$ is denoted by

$$\omega_k(\rho, x, \delta)_{p,\beta} = \sup\{|\Delta_r^k \rho(x)| : x, x + rk \in [x - \tfrac{\delta k}{2}, x + \tfrac{\delta k}{2}]\} \qquad (3)$$

such that

$$\Delta_r^k(\rho, x) = \sum_{r=0}^k (-1)^{r+k} \binom{k}{r} \rho(x + r\delta) \text{ if } x, x + r\delta \in X \qquad (4)$$

The average modulus of smoothness we defined by

$$\tau_k(\xi, \delta)_{p,\beta} = \|\omega_k(\rho, ., \delta)\|_{p,\beta}. \qquad (5)$$

Let $\mathbb{N}$ be the set of natural numbers and $\mathbb{H}_k$ the set of all algebraic polynomials of degree less than or equal to $k \in \mathbb{N}$.

The degree of best one-sided approximation in $L_{p,\beta}(X), 1 \leq p < \infty$, of unbounded function as $\rho$ by operators $\mathcal{P}_k \& q_k$ are denoted by:

$$\tilde{\mathcal{E}}_k(\rho, x)_{p,\beta} = \inf\{\|\rho - \mathcal{P}_k\|_{p,\beta} : \mathcal{P}_k \in \mathbb{H}_k\} \qquad (6)$$

(٨٦)

$$\tilde{\mathcal{E}}_k(\rho,x)_{p,\beta} = \inf\{\|q_k - \mathcal{P}_k\|_{p,\beta} : \mathcal{P}_k, q_k \in \mathbb{H}_k \,\&\, \mathcal{P}_k(x) \leq \rho(x) \leq q_k(x)\}. \quad (7)$$

It easy to verify that there are not linear operators for one-sided approximation in X. Some non-linear construction have been proposed in [ 3] and [ 6 ].

Let us consider the step function

$$\Phi(x) = \begin{cases} 0, & if -1 \leq x \leq 0, \\ 1, & if \ 0 < x \leq 1, \end{cases} \quad (8)$$

determine two collections of functions $\{\mathcal{P}_k\}$ and $\{q_k\}$, $\mathcal{P}_k, q_k \in \mathbb{H}_k$ such that

$$\mathcal{P}_k(x) \leq \Phi(x) \leq q_k(x), \ x \in [-1,1] \quad (9)$$

and

$$C_k = \|\rho - \mathcal{P}_k\|_{p,\beta} \to 0, \ p = 1 \quad (10)$$

For the first one we work in space $L_{p,\beta}(X)$. For $1 \leq p < \infty$, we construct two different sequences of operators, for $x \in X$, $k \in \mathbb{N}$ and $\rho \in L_{p,\beta}(X)$ define

$$\mathcal{M}_k(\rho,x) = \rho(0) + \int_0^1 \mathcal{P}_k(t-x)\, \rho'_+(t)dt - \int_0^1 q_k(t-x)\, \rho'_-(t)dt \quad (11)$$

and

$$\mathcal{N}_k(\rho,x) = \rho(0) + \int_0^1 q_k(t-x)\, \rho'_+(t)dt - \int_0^1 \mathcal{P}_k(t-x)\, \rho'_-(t)dt \quad (12)$$

it is clear $\mathcal{M}_k(\rho), \mathcal{N}_k(\rho) \in \mathbb{H}_k$, we will prove that

$$\mathcal{M}_k(\rho,x) \leq \rho(x) \leq \mathcal{N}_k(\rho,x), x \in X \text{ and both}$$

$$\|\rho - \mathcal{M}_k(\rho)\|_{p,\beta} \leq C_k \|\rho'\|_{p,\beta} \text{ and}$$

$$\|\rho - \mathcal{N}_k(\rho)\|_{p,\beta} \leq C_k \|\rho'\|_{p,\beta}, \text{ where } C_k \text{ be as in (10)}.$$

In the second part, for function $L_{p,\beta}(X)$, we construct operators

$$G_y(\rho,x) = \int_0^1 [\rho(1-y)x + yt) - \omega(\rho,(1-y)x + yt, y)]dt \quad (13)$$

and

$$H_y(\rho,x) = \int_0^1 [\rho(1-y)x + yt) + \omega(\rho,(1-y)x + yt, y)]dt. \quad (14)$$

It is clear $G_k(\rho), H_k(\rho) \in \mathbb{H}_k$ and so we can define

$$\mathcal{L}_{k,y}(\rho,x) = \mathcal{M}_k(G_y(\rho),x) \quad (15)$$

$$\mathcal{J}_{k,y}(\rho,x) = \mathcal{N}_k(H_y(\rho),x), \qquad (16)$$

Where $\mathcal{M}_k(\rho)$ and $\mathcal{N}_k(\rho)$ by as equations (11) and (12) in that order.

We will prove that

$$\mathcal{L}_{k,y}(\rho,x) \leq \rho(x) \leq \mathcal{J}_{k,y}(\rho,x), \ x \in X \text{ and}$$

present the degree of best one-sided approximation of unbounded functions by operators $\mathcal{L}_{k,y}(\rho,x)$ and $\mathcal{J}_{k,y}(\rho,x), x \in X$ in terms averaged modulus of continuity.

In the last years there has been interest in studying open problems related to one-sided approximations (see [1] , [2]).

We point out that other operators for one-sided approximations have constructed in [7].

In particular, the operators presented in [6] yield the non-optimal rate $O(\tau\left(\rho,\frac{1}{\sqrt{k}}\right))$ where is ones consider in [4] give the optimal rate, but without an explicit constant. The paper is organized as follows. In section (3) we calculate the degree of best one-sided approximation of unbounded functions by mean of the operators define (13) and (14).

Finally in the some section, we consider the degree of the best onesided approximation by mean of the operators defined in (15) and (16)

### 3.AUXILIARY LEMMAS

**Lemma 3.1:**

Let $\rho \in L_{p,\beta}[0,1], \ y \in (0,1)$, and the operators $G_y(\rho) \ \& \ H_y(\rho)$ are defined through equations (13) and (14) in that order.

Then, $G_y(\rho,x) \leq \rho(x) \leq H_y(\rho,x), \ x \in X = [0,1]$ and

$$Max\{\|G'_y\|_{p,\beta}, \|H'_y\|_{p,\beta}\} \leq \frac{3}{k}\tau_k(\rho,y)_{p,\beta}.$$

**Lemma 3.2:**

Let $\Phi(x)$ be given in equation (8). Every $x \in [-1,1]$ we known

$$\mathcal{P}_k(\rho) = T_k^-(arc \ cosx) \text{ and}$$

$$q_k(\rho) = T_k^+(arc \ cosx).$$

Then, $\mathcal{P}_k, q_k \in \mathbb{H}_k, \mathcal{P}_k(\rho) \leq \rho(x) \leq q_k(\rho), \ x \in [-1,1]$ and

$$\|q_k(\rho) - \mathcal{P}_k(\rho)\|_{p,\beta} \leq \frac{4\pi^2}{k+2}.$$

**Lemma 3.3:**

Let $\rho \in L_{p,\beta}(X), 1 \leq p < \infty$, $k \in \mathbb{N}$, and $k \geq 2$. Let $\mathcal{M}_k(\rho)$ and $\mathcal{N}_k(\rho)$ by as equations (11) and (12) in that order. Then,

$\mathcal{M}_k(\rho), \mathcal{N}_k(\rho) \in \mathbb{H}_k$ and

$\mathcal{M}_k(\rho, x) \leq \rho(x) \leq \mathcal{N}_k(\rho, x), x \in X = [0,1]$.

*Proof:*

From equations (9), (10), (11) and (12), it is clear that

$\mathcal{M}_k(\rho), \mathcal{N}_k(\rho) \in \mathbb{H}_k$. We have

$\mathcal{M}_k(\rho, x) = \rho(0) + \int_0^1 \mathcal{P}_k(t - x)\rho'_+(t)dt - \int_0^1 q_k(t - x)\rho'_-(t)dt$

where $\mathcal{P}_k, q_k \in \mathbb{H}_k$, such that $\mathcal{P}_k(\rho) \leq \rho(x) \leq q_k(\rho)$, $x \in X$

and $\|\mathcal{P}_k - q_k\|_{p,\beta} \to 0$

since, $\mathcal{P}_k(\rho) \leq \rho(x) \leq q_k(\rho), x \in [0,1]$,

thus

$\mathcal{M}_k(\rho, x) \leq \rho(0) + \int_0^1 \Phi(t - x)\rho'_+(t)dt - \int_0^1 \Phi(t - x)\rho'_-(t)dt$

$= \rho(0) + \int_0^1 \Phi(t - x)\rho'(t)dt$

$= \rho(0) + \rho(x) - \rho(0)$

$= \rho(x)$.

So,

$\rho(x) = \rho(0) + \rho(x) - \rho(0)$

$= \rho(0) + \int_0^1 \rho'(t)dt$

$= \rho(0) + \int_0^1 \Phi(t - x)\rho'(t)dt$

$= \rho(0) + \int_0^1 \Phi(t - x)\rho'_+(t)dt - \int_0^1 \Phi(t - x)\rho'_-(t)dt$

$\leq \rho(0) + \int_0^1 q_k(t - x)\rho'_+(t)dt - \int_0^1 \mathcal{P}_k(t - x)\rho'_-(t)dt$

$= \mathcal{N}_k(\rho, x)$.

**Lemma 3.4:**

For $\rho \in L_{p,\beta}(X), 1 \leq p < \infty$, $k \in \mathbb{N}$, and $k \geq 2$. Let $\mathcal{M}_k(\rho)$ and $\mathcal{N}_k(\rho)$ by as equations (11) and (12) in that order. Then,

$$Max\{\|\rho - \mathcal{M}_k(\rho)\|_{p,\beta}, \|\rho - \mathcal{N}_k(\rho)\|_{p,\beta}\} \leq C_k \|\rho'\|_{p,\beta}.$$

*Proof:*

Since

$$|\rho(x) - \mathcal{M}_k(\rho, x)| \leq \int_{-x}^{1-x}(q_k(h) - \mathcal{P}_k(h))|\rho'(x+h)|dh,$$

We put $\alpha_k(h) = q_k(h) - \mathcal{P}_k(h)$ and from Holder's inequity, we get

$$(\|\rho - \mathcal{M}_k(\rho)\|_{p,\beta})^p \leq \int_0^1 \left|\frac{\int_{-x}^{1-x} \alpha_k(h)|\rho'(x+h)|dh}{\beta(x)}\right|^p dx$$

$$\leq \int_0^1 \left(\left|\int_{-x}^{1-x} \alpha_k(h)dh\right|^{p-1}\right)\left(\left|\frac{\int_{-x}^{1-x} \alpha_k(h)|\rho'(x+h)|^p dh}{\beta(x)}\right|\right) dx$$

$$\leq \left(\int_0^1 |\alpha_k(u)|^{p-1} du\right)\left(\int_0^1 \left|\frac{\rho'(v)}{\beta(v)}\right|^p \left(\int_{v-1}^v \frac{\alpha_k(h)}{\beta(h)} dh\right) dv\right)$$

$$\leq \left(\int_0^1 |\alpha_k(u)|^p du\right)\left(\int_0^1 \left|\frac{\rho'(v)}{\beta(v)}\right|^p dv\right)$$

Thus

$$\|\rho - \mathcal{M}_k(\rho)\|_{p,\beta} \leq \left(\int_0^1 |\alpha_k(u)|^p du\right)^{1/p}\left(\int_0^1 \left|\frac{\rho'(v)}{\beta(v)}\right|^p dv\right)^{1/p},$$

Hence

$$\|\rho - \mathcal{M}_k(\rho)\|_{p,\beta} \leq \|\alpha_k\|_p \|\rho'\|_{p,\beta} = C_k \|\rho'\|_{p,\beta}.$$

Likewise, we show that,

$$\|\rho - \mathcal{N}_k(\rho)\|_{p,\beta} \leq C_k \|\rho'\|_{p,\beta}.$$

### 4. MAIN RESULTS

**Theorem 4.1:**

Let $\rho \in L_{p,\beta}(X), 1 \leq p < \infty,\ k \in \mathbb{N}$, and $k \geq 2$. Let $G_y(\rho)\ \&\ H_y(\rho)$ are defined through equations (13) and (14) in that order.

Then,

$$Max\{\|\rho - G_y(\rho)\|_{p,\beta}, \|\rho - H_y(\rho)\|_{p,\beta}\} \leq C_1(y,p)\tau_k(\rho,y)_{p,\beta}$$

and

$$\tilde{\mathcal{E}}_k(\rho)_{p,\beta} \leq C_k(y,p)\tau_k(\rho,y)_{p,\beta}.$$

*Proof:*

It is used to, take $q$ such that $1/p + 1/q = 1$, from equations (13), (14) and Holder's inequality, we get

$$(y\|\rho - G_y(\rho)\|_{p,\beta})^p = y^p \int_0^1 \left|\frac{\rho(x) - G_y(\rho,x)}{\beta(x)}\right|^p dx$$

$$\leq y^p \int_0^1 \left|\frac{H_y(\rho,x) - G_y(\rho,x)}{\beta(x)}\right|^p dx$$

$$\leq 2^p y^p \int_0^1 \int_0^y \left|\frac{\omega(\rho,(1-y)x+yt,y)}{\beta((1-y)x)}\right|^p dt dx.$$

Put $h = (1-y)x$ implies $dh = (1-y)dx$

$$(y\|\rho - G_y(\rho)\|_{p,\beta})^p \leq \frac{2^p y^p}{1-y} \int_0^y \int_t^{1-y+t} \left|\frac{\omega(\rho,h,y)}{\beta(h)}\right|^p dh dt$$

$$\leq \frac{2^p y^{p/q}}{1-y} \int_0^y \int_0^1 \left|\frac{\omega(\rho,h,y)}{\beta(h)}\right|^p dh dt$$

$$\leq \frac{2^p y^{\frac{p}{q}+1}}{1-y} \int_0^1 \left|\frac{\omega(\rho,h,y)}{\beta(h)}\right|^p dh$$

Thus,

$$\|\rho - G_y(\rho)\|_{p,\beta} \leq \frac{2}{(1-y)^{1/p}} \left(\int_0^1 \left|\frac{\omega(\rho,h,y)}{\beta(h)}\right|^p dh\right)^{1/p}$$

$$\leq \frac{2}{(1-y)^{1/p}} \|\omega(\rho,.,y)\|_{p,\beta}$$

$$= \frac{2}{(1-y)^{1/p}} \tau_k(\rho,y)_{p,\beta}$$

We have $\frac{2}{(1-y)^{1/p}}$ constant depending on $y$ and $p$, then

$$\|\rho - G_y(\rho)\|_{p,\beta} \leq C_1(y,p)\tau_k(\rho,y)_{p,\beta}.$$

Analogously, we can show $\|\rho - H_y(\rho)\|_{p,\beta} \leq C_1(y,p)\tau_k(\rho,y)_{p,\beta}.$

We go to the following inequality:

$$\tilde{\mathcal{E}}_k(\rho)_{p,\beta} \leq \|H_y(\rho) - G_y(\rho)\|_{p,\beta}$$

$$\leq \|\rho - H_y(\rho)\|_{p,\beta} + \|\rho - G_y(\rho)\|_{p,\beta}$$

$$\leq C_k(y,p)\tau_k(\rho,y)_{p,\beta}.$$

**Theorem 4.2:**

Let $\rho \in L_{p,\beta}(X), 1 \leq p < \infty, k \in \mathbb{N}$. Let $\mathcal{L}_{k,y}(\rho)$ and $\mathcal{J}_{k,y}(\rho)$ are defined through equations (13) and (14) in that order.

Then,

$$\mathcal{L}_{k,y}(\rho,x) \le \rho(x) \le \mathcal{J}_{k,y}(\rho,x), \ x \in X$$

$$Max\{\|\rho - \mathcal{L}_{k,y}(\rho)\|_{p,\beta}, \|\rho - \mathcal{J}_{k,y}(\rho)\|_{p,\beta}\} \le (C_1(y,p) + \frac{3C_k}{y})\tau_k(\rho,y)_{p,\beta}$$

and

$$\tilde{\mathcal{E}}_k(\rho)_{p,\beta} \le (C_1(y,p) + \frac{6C_k}{y})\tau_k(\rho,y)_{p,\beta}.$$

**_Proof:_**

Let $G_y(\rho)$ and $H_y(\rho)$ are defined through equations (13) and (14) in that order.

So, from equations (15) and (16), it is clear $\mathcal{L}_{k,y}(\rho), \mathcal{J}_{k,y}(\rho) \in \mathbb{H}_k$.

Also, from equations (15), (16), Theorem 4.1, Lemma 3.3, Lemma 3.4 and Lemma 3.1, since

$$\mathcal{L}_{k,y}(\rho,x) = \mathcal{M}_k\left(G_y(\rho),x)\right) \le G_y(\rho,x) \le \rho(x)$$

$$\le H_y(\rho,x) \le \mathcal{N}_k\left(H_y(\rho),x)\right) = \mathcal{J}_{k,y}(\rho,x), \ x \in [0,1].$$

Moreover,

$$\|\rho - \mathcal{L}_{k,y}(\rho)\|_{p,\beta} \le \|\rho - G_y(\rho)\|_{p,\beta} + \|G_y(\rho) - \mathcal{L}_{k,y}(\rho)\|_{p,\beta}$$

$$\le C_1(y,p)\tau_k(\rho,y)_{p,\beta} + \|\rho - \mathcal{M}_k\left(G_y(\rho,x)\right)\|_{p,\beta}$$

$$\le C_1(y,p)\tau_k(\rho,y)_{p,\beta} + C_k\|G'_y(\rho)\|_{p,\beta}$$

$$= C_1(y,p)\tau_k(\rho,y)_{p,\beta} + C_k\left\|\frac{G'_y(\rho,.)}{\beta(.)}\right\|_p$$

$$\le C_1(y,p)\tau_k(\rho,y)_{p,\beta} + \frac{3C_k}{y}\tau_k(\frac{\rho}{\beta},y)_p$$

$$= C_1(y,p)\tau_k(\rho,y)_{p,\beta} + \frac{3C_k}{y}\tau_k(\rho,y)_{p,\beta}$$

$$= (C_1(y,p) + \frac{3C_k}{y})\tau_k(\rho,y)_{p,\beta}.$$

The approximation for $\|\rho - \mathcal{J}_{k,y}(\rho,x)\|_{p,\beta}$ follows similarly.

Thus,

$$\tilde{\mathcal{E}}_k(\rho)_{p,\beta} \le \|\mathcal{J}_{k,y}(\rho,x) - \mathcal{L}_{k,y}(\rho)\|_{p,\beta}$$

$$\le \|\rho - \mathcal{L}_{k,y}(\rho)\|_{p,\beta} + \|\rho - \mathcal{J}_{k,y}(\rho,x)\|_{p,\beta}$$

$$\leq 2(C_1(y,p) + \frac{3C_k}{y})\tau_k(\rho,y)_{p,\beta}$$

$$\leq (C_1(y,p) + \frac{6C_k}{y})\tau_k(\rho,y)_{p,\beta}.$$

**Theorem 4.3:**

Let $\rho \in L_{p,\beta}(X), 1 \leq p < \infty, \ k \in \mathbb{N},$ and $k \geq 2.$ Let $\mathcal{P}_k(\rho)$ and $q_k(\rho)$ be the sequence of polynomials constructed as in (9), set

$$A_k(\rho) = \mathcal{L}_{k,\frac{1}{k}}(\rho) \text{ and } B_k(\rho) = \mathcal{J}_{k,\frac{1}{k}}(\rho),$$

where

$$\mathcal{L}_{k,\frac{1}{k}}(\rho) \text{ and } \mathcal{J}_{k,\frac{1}{k}}(\rho) \text{ are given as (15) and (16) respectively.}$$

Then

$$A_k(\rho,x) \leq \rho(x) \leq B_k(\rho,x), \ x \in X,$$

$$Max \ \{\|\rho - A_k(\rho)\|_{p,\beta}, \|\rho - B_k(\rho)\|_{p,\beta}\} \leq (C_1(y,p) + \frac{3C_k}{y})\tau_k(\rho,\frac{1}{k})_{p,\beta}$$

and

$$\tilde{\mathcal{E}}_k(\rho)_{p,\beta} \leq 2(C_k(y,p) + \frac{12k\pi^2}{k+2})\tau_k(\rho,\frac{1}{k})_{p,\beta}.$$

***Proof:***

Using equations (15) and (16) with $y = \frac{1}{k}$ and $k \geq 2,$ we get

$$\mathcal{L}_{k,\frac{1}{k}}(\rho,x) = \mathcal{M}_k\left(G_{\frac{1}{k}}(\rho,x)\right) \text{ and } \mathcal{J}_{k,\frac{1}{k}}(\rho,x) = \mathcal{N}_k\left(H_{\frac{1}{k}}(\rho,x)\right)$$

where

$$G_{\frac{1}{k}}(\rho), H_{\frac{1}{k}}(\rho) \in \mathbb{H}_k. \text{ Also}$$

$$\mathcal{M}_k\left(G_{\frac{1}{k}}(\rho)\right), \mathcal{N}_k\left(H_{\frac{1}{k}}(\rho)\right) \in \mathbb{H}_k$$

Using Lemma 3.3, since $\mathcal{M}_k(\rho,x) \leq \rho(x) \leq \mathcal{N}_k(\rho,x), \ x \in X.$

Hence, $A_k(\rho,x) \leq \rho(x) \leq B_k(\rho,x), \ x \in X.$

We need an approximate for $\|\rho - A_k(\rho)\|_{p,\beta}$ one has:

Using (15), Lemma 3.2 and Theorem 4.2 we obtain

$$\|\rho - A_k(\rho)\|_{p,\beta} = \left\|\rho - \mathcal{L}_{k,\frac{1}{k}}(\rho)\right\|_{p,\beta} \leq (C_k(y,p) + \frac{3C_k}{\frac{1}{k}})\tau_k(\rho,\frac{1}{k})_{p,\beta}$$

$$\leq (C_k(y,p) + \frac{12k\pi^2}{k+2})\tau_k(\rho,\frac{1}{k})_{p,\beta}.$$

Analogously, we can show

$$\|\rho - B_k(\rho)\|_{p,\beta} \leq (C_k(y,p) + \frac{12k\pi^2}{k+2})\tau_k(\rho,\frac{1}{k})_{p,\beta}.$$

Thus,

$$\tilde{\mathcal{E}}_k(\rho)_{p,\beta} \leq \|B_k(\rho) - A_k(\rho)\|_{p,\beta}$$

$$\leq \|B_k(\rho) - \rho\|_{p,\beta} + \|\rho - A_k(\rho)\|_{p,\beta}$$

$$\leq 2(C_k(y,p) + \frac{12k\pi^2}{k+2})\tau_k(\rho,\frac{1}{k})_{p,\beta}.$$